\RequirePackage{pdf14}

\documentclass[preprint,3p,12pt]{elsarticle}

\usepackage[english]{babel}
\usepackage{amsmath}
\usepackage{amsthm}
\usepackage{amstext}
\usepackage{amssymb}
\usepackage{mathtools}

\usepackage{varioref}
\usepackage{hyperref}

\usepackage[noabbrev,capitalise,nameinlink]{cleveref}

\hypersetup{colorlinks={true},linkcolor={blue},citecolor=green}

\usepackage{amssymb,nicefrac,bm,upgreek,mathtools,verbatim}
\usepackage[mathscr]{eucal}
\usepackage{dsfont}
\RequirePackage[normalem]{ulem}
\usepackage{tikz}
\usepackage{color}

\allowdisplaybreaks

\newtheorem{theorem}{Theorem}[section]

\newtheorem{proposition}{Proposition}[section]
\newtheorem{lemma}{Lemma}[section]

\theoremstyle{definition}
\newtheorem{definition}{Definition}[section]

\newtheorem{assumption}{Assumption}[section]

\newcommand{\stkout}[1]{\ifmmode\text{\sout{\ensuremath{#1}}}\else\sout{#1}\fi}

\newcommand{\ttup}[1]{\textup{(}#1\textup{)}}

\newcommand{\smid}{\,|\,}               
\newcommand{\df}{\coloneqq}             
\newcommand{\Ind}{\bm{1}}               
\DeclareMathOperator{\Exp}{\mathbb{E}}  
\DeclareMathOperator{\Prob}{\mathbb{P}} 
\newcommand{\D}{\mathrm{d}}             
\newcommand{\RR}{\mathbb{R}}            
\newcommand{\NN}{\mathbb{N}}            
\DeclareMathOperator*{\supp}{support}

\newcommand{\Act}{\mathbb{U}}           
\newcommand{\Uadm}{\Gamma_{\mathsf{A}}} 
\newcommand{\Um}{\Gamma_{\mathsf{M}}}   
\newcommand{\Usm}{\Gamma_{\mathsf{S}}}  
\newcommand{\Usd}{\Gamma_{\mathsf{SD}}} 

\newcommand{\cG}{\mathcal{G}}           
\newcommand{\cH}{\mathcal{H}}           
\newcommand{\cT}{\mathcal{T}}           
\newcommand{\Bor}{{\mathcal{B}}}        
\newcommand{\Pm}{{\mathcal{P}}}         
\newcommand{\cM}{{\mathcal{M}}}         

\newcommand{\XX}{{\mathbb{X}}}          
\newcommand{\KK}{\mathbb{K}}            
\newcommand{\HH}{\mathbb{H}}            
\newcommand{\cU}{{\mathcal{U}}}         

\newcommand{\abs}[1]{\lvert#1\rvert}

\newcommand{\babs}[1]{\bigl\lvert#1\bigr\rvert}

\newcommand{\Babss}[1]{\Biggl\lvert#1\Biggr\rvert}
\newcommand{\bnorm}[1]{\bigl\lVert#1\bigr\rVert}

\definecolor{dmagenta}{rgb}{.4,.1,.5}
\definecolor{dblue}{rgb}{.0,.0,.5}
\definecolor{mblue}{rgb}{.0,.0,.7}
\definecolor{ddblue}{rgb}{.0,.0,.4}
\definecolor{dred}{rgb}{.7,.0,.0}
\definecolor{dgreen}{rgb}{.0,.5,.0}
\definecolor{Eeom}{rgb}{.0,.0,.5}

\tikzstyle{terminal}=[rectangle, rounded corners, minimum width=1cm, minimum height=1cm,text centered, draw=black]

\tikzstyle{terminal_2}=[rectangle, rounded corners, minimum width=0.5cm, minimum height=0.5cm,text centered, draw=black]

\tikzstyle{terminal_3}=[ellipse, rounded corners, minimum width=0.5cm, minimum height=0.5cm,text centered, draw=black, ultra thick]

\tikzstyle{arrow} = [thick,->,>=stealth]

\tikzstyle{addblock} = [draw,circle]

\tikzstyle{cost_block}=[rectangle, minimum width=1cm, minimum height=1cm,text centered, draw=black, ultra thick]

\allowdisplaybreaks

\begin{document}

\begin{frontmatter}

\title{Convex Analytic Method Revisited: Further Optimality Results and Performance of Deterministic Policies in Average Cost Stochastic Control \tnoteref{label0}}
\tnotetext[label0]{This research was partially supported by the Natural Sciences and Engineering Research Council of Canada (NSERC), by the National Science Foundation through grant DMS-1715210,
and by the Army Research Office through grant W911NF-17-1-001.}

\author[label1]{Ari Arapostathis}
\address[label1]{University of Texas at Austin.}
\ead{ari@utexas.edu}

\author[label2]{Serdar Y\"uksel\corref{cor}}
\ead{yuksel@queensu.ca}
\cortext[cor]{Corresponding author}
\address[label2]{Department of Mathematics and Statistics, Queen's University, Canada.}

\begin{abstract}
The convex analytic method has proved to be a very versatile method
for the study of infinite horizon average cost optimal stochastic control problems.
In this paper, we revisit the convex analytic method and make three primary contributions:
(i) We present an existence result for controlled Markov models that lack weak continuity of the transition kernel but are strongly continuous in the action variable for every fixed state variable.
(ii) For average cost stochastic control problems in standard Borel spaces,
while existing results establish the optimality of stationary (possibly randomized) policies,
few results are available on the optimality of deterministic policies. We review existing results and present further conditions under which an average cost optimal stochastic
control problem admits optimal solutions that are deterministic stationary. (iii) We establish conditions under which the performance under stationary deterministic (and also quantized) policies is dense in the set of performance values under randomized stationary policies.

\end{abstract}


\end{frontmatter}

\section{Introduction}



We start by reviewing the usual model in the literature for controlled
Markov chains, otherwise referred to as Markov decision processes (MDPs).
In general, for a topological space $\mathcal{X}$, we denote by $\Bor(\mathcal{X})$
its Borel $\sigma$-field and by $\Pm(\mathcal{X})$ the set of probability
measures on $\Bor(\mathcal{X})$.

A \emph{controlled Markov chain} consists of the tuple $\bigl(\XX,\Act,\cU,\cT,c\bigr)$,
whose elements can be described as follows.
\begin{itemize}
\item[(a)]
The \emph{state space} $\XX$ and the
\emph{action} or
\emph{control space} $\Act$ are Borel subsets of complete, separable,
metric (i.e., \emph{Polish}) spaces.
\item[(b)]
The map $\cU\colon \XX\to \Bor(\Act)$ is a strict, measurable
multifunction.
The set of admissible state/action pairs is 
\begin{equation*}
\KK \,\df\, \bigl\{(x,u)\colon\, x\in\XX,\,u\in\cU(x)\bigr\}\,,
\end{equation*}
endowed with the subspace topology
corresponding to $\Bor(\XX\times\Act)$.
\item[(c)]
The map $\cT\colon\KK\to\Pm(\XX)$ is a stochastic kernel on
$\KK\times\Bor(\XX)$, that is, $\cT(\,\cdot\smid  x,u)$ is
a probability measure on $\Bor(\XX)$ for each $(x,u)\in\KK$,
and $(x,u) \mapsto \cT(A\smid  x,u)$ is measurable for each $A\in\Bor(\XX)$.
\item[(d)]
The map $c\colon \KK\to \mathbb{R}_+$ is measurable, and is called the \emph{running
cost} or \emph{one stage cost}.
We assume that it is bounded from below in $\KK$, so without loss of generality,
it takes values in $[1,\infty)$.
\end{itemize}

The (admissible) \emph{history spaces} are defined as
\begin{equation*}
\HH_0\,\df\, \XX\,, \quad \HH_{t} \,\df\, \HH_{t-1}\times \mathbb{U} \times \mathbb{X} ,\quad t\in\NN\,,
\end{equation*}
and the canonical sample space is defined as $\Omega\df (\XX\times\Act)^\infty$.
These spaces are endowed with their respective product topologies and
are therefore Borel spaces.
The state, action (or control), and information processes, denoted
by $\{X_t\}_{t\in\NN_0}$,  $\{U_t\}_{t\in\NN_0}$ and $\{H_t\}_{t\in\NN_0}$,
respectively, are defined by the projections
\begin{equation*}
X_t(\omega) \,\df\, x_t\,,\quad U_t(\omega) \,\df\, u_t\,,
\quad H_t(\omega) \,\df\, (x_0,u_0, \dotsc, u_{t-1}, x_t)
\end{equation*}
for each $\omega=(x_0,u_0, \dotsc,u_{t-1},x_t, u_t,\dotsc)\in\Omega$.
An \emph{admissible control policy}, or \emph{policy}, is a sequence
$\gamma = \{\gamma_t\}_{t\in\NN_0}$ of stochastic kernels on
$\HH_t\times\Bor(\Act)$ satisfying the constraint
\begin{equation*}
\gamma_t(\cU(X_t)\mid h_t) \,=\, 1\,,\quad  h_t\in\HH_t\,.
\end{equation*}
The set of all admissible policies is denoted by $\Uadm$.
It is well known
(see \cite[Prop.\ V.1.1, pp.~162--164]{Neveu})
that for any given $\nu\in\Pm(\XX)$ and $\gamma\in\Uadm$ there exists
a unique probability measure $\Prob^\gamma_\nu$ on $\bigl(\Omega,\Bor(\Omega)\bigr)$
satisfying
\begin{align*}
\Prob^\gamma_\nu(X_0\in D) &\,=\, \nu(D)\qquad\forall\, D\in\Bor(\XX)\,,\\
\Prob^\gamma_\nu(U_t\in C\mid H_t) &\,=\, 
\gamma_t(C\mid H_t)
\quad \Prob^\gamma_\nu\text{-a.s.}\,,\quad \forall\, C\in\Bor(\Act)\\
\Prob^\gamma_\nu(X_{t+1}\in D\mid H_t,U_t)  &\,=\,
\cT(D\mid X_t, U_t) \quad \Prob^\gamma_\nu\text{-a.s.}\,,
\quad\forall\, D\in\Bor(\XX)\,.
\end{align*}
The expectation operator corresponding to $\Prob^\gamma_\nu$
is denoted by $\Exp^\gamma_\nu$.
If $\nu$ is a Dirac mass at $x\in\XX$, we simply write these
as $\Prob^\gamma_x$ and $\Exp^\gamma_x$.

A policy $\gamma$ is called \emph{Markov}
if there exists a sequence of measurable maps $\{v_t\}_{t\in\NN_0}$,
where $v_t\colon \XX\to \Pm(\Act)$, where $\Pm(\Act)$ is endowed with the weak convergence topology, for each $t\in\NN_0$,
such that
$$\gamma_t(\,\cdot\mid H_t) \,=\, v_t(X_t)(\cdot) \quad \Prob^\gamma_\nu\text{-a.s.}$$
With some abuse in notation, such a policy is identified
with the sequence $v=\{v_t\}_{t\in\NN_0}$.
Note then that
 $\gamma_t$ may be written as a stochastic kernel $\gamma_t(\cdot\smid  x)$ on
$\XX\times\Bor(\Act)$ which satisfies $\gamma_t(\cU(x)\smid  x)=1$.
Let $\Um$ denote the set of all Markov policies. 

We say that a  Markov policy $\gamma$ is \emph{deterministic}, or \emph{simple},
if $\gamma_t$ is a Dirac mass, in which case $\gamma_t$ is identified
with a Borel measurable function $\gamma_t\colon \XX\to\Act$.
In other words, $\gamma_t$ is a measurable selector from the set-valued map
$\cU(x)$ \cite{Feinberg-13}. We let $\Gamma_{\mathsf{MD}}$ denote the set of deterministic Markov policies.

We add the adjective \emph{stationary} to indicate that the Markov policy
does not depend on $t\in\NN_0$, that is, $\gamma_t=\gamma$ for all $t\in\NN_0$.
We let $\Usm$ denote the class of stationary Markov policies,
henceforth referred to simply as \emph{stationary policies},
and let $\Usd\subset\Usm$ denote the subset of those that are deterministic.

In summary, under a policy $\gamma\in\Usm$, the process
satisfies the following: for all Borel sets $B \in \mathcal{B}(\XX)$, $t \ge 0$, and
($\Prob^\gamma$ almost all) realizations $X_{[0,t]}, U_{[0,t]}$, we have
\begin{equation} \label{eq_evol}
\begin{aligned}
\Prob^\gamma\bigl( X_{t+1} \in B \smid  X_{[0,t]}=x_{[0,t]}, U_{[0,t]}=u_{[0,t]}\bigr)
&\,=\, \Prob^\gamma( X_{t+1} \in B \smid  X_t=x_t, U_t=u_t) \\
&\,= \cT(  B \smid  x_t, u_t)\,.
\end{aligned}
\end{equation}

Using stochastic realization results (see  \cite[Lemma~1.2]{gihman2012controlled},
or \cite[Lemma~3.1]{BorkarRealization}), stochastic processes that satisfy \cref{eq_evol} admit a realization in the form 
\begin{equation}\label{EE1.2}
X_{t+1}\,=\,f(X_t,U_t,W_t)
\end{equation}
almost surely,
where $f$ is measurable and $w_t$ is i.i.d. $[0,1]-$valued.
Since a system of the form \cref{EE1.2} satisfies \cref{eq_evol},
it follows that the representations in these equations are equivalent.

In this paper, we consider the problem of minimizing the average cost
\begin{equation}\label{E-AC}
J^*(x)\,\df\, \inf_{\gamma \in \Uadm} J(x,\gamma) \,=\, \inf_{\gamma \in \Uadm}\,\limsup_{T \to \infty}\,
\frac{1}{T}\, \Exp^{\gamma}_x \Biggl[\sum_{t=0}^{T-1} c(X_t,U_t)\Biggr]\,.
\end{equation}
We say that a policy $\gamma\in\Uadm$ is optimal if it attains the infimum in \cref{E-AC}.

This is an important problem in applications where one is concerned about the long-term behaviour, unlike the discounted cost setup where the primary interest is in the short-term time stages. 
 
For the study of the average cost problem, there are three commonly adopted approaches \cite{survey}; contraction or value iteration based methods (see e.g.
\cite{Veg03,hernandez2012adaptive,ABor-19}), the vanishing discount method (see e.g. \cite{survey,feinberg2012average,HernandezLermaMCP,hernandezlasserre1999further,costa2012average,GoHe95,yu2020average} which have various conditions and relaxations),
and the convex analytic method (to be reviewed further below).
The first two are based on the arrival at what what is known as the
\emph{average cost optimality equation} (ACOE)
(and its variation involving an inequality (ACOI)).
Efforts under this method typically (and as we will study, not necessarily) require some
recurrence/ergodicity/Dobrushin type geometric or at least subgeometric convergence conditions,
which may be too strong for a large class of applications
(e.g., for belief-MDP reduction of Partially Observable Markov Decision Processes). 

The third approach, via the convex analytic method, is based on the properties of
expected (or sample path) occupation measures and their limit behaviours,
leading to a linear program involving the space of probability measures. The convex analytic approach, typically attributed to Manne \cite{Manne} and Borkar \cite{Borkar2} (see additionally \cite{kurano1989existence,hernandez1993existence,survey,HernandezLermaMCP,yu2020minimum}), is a versatile approach to the optimization of infinite-horizon problems, which leads to a linear program. This approach is particularly effective for constrained optimization problems and infinite horizon average cost optimization problems. It avoids the use of dynamic programming and can also be tailored towards obtaining results on sample-path optimality via martingale convergence theorems under mild continuity conditions \cite{vega1999sample,lasserre1999sample,arapostathis2017some,yu2020minimum}. Most importantly perhaps, this approach generally requires less restrictive conditions on the existence of an optimal policy for average cost stochastic control. 

These approaches are related through a duality analysis, as noted in
\cite[Chapter 6]{HernandezLermaMCP} (see also \cite{hernandez1993existence} for a direct argument under positive Harris recurrence assumptions). 
However the more general conditions leading to solutions under these approaches are not identical,
therefore, the corresponding conditions of existence and structural results for optimal
policies are somewhat different. That is, going from one approach to another one
(e.g., from the convex analytic solution to an ACOE) still entails open problems. 

For MDPs with weakly or strongly continuous transition kernels, if ACOE/ACOI can be established (under somewhat strong conditions as reviewed above), the existence of deterministic stationary optimal policies naturally follows. While the convex analytic method typically provides less conservative conditions for existence of optimal policies, whether the optimal policy can be taken to be {\it deterministic} is generally an open question with only few results reported in the literature. This question is a further primary motivation for this paper. Optimality of deterministic policies finds itself in many applications, e.g. in the optimal zero-delay quantization problem \cite{BorkarMitterTatikonda,YukLinZeroDelay} for average cost criterion, where common randomness between an encoder and decoder would be costly to implement. 

\medskip\noindent{\bf Contributions.} 
\begin{itemize}
\item[(i)] In Theorem \ref{thmSetwise}, we present an existence result for average cost controlled Markov models that are strongly continuous in the action for every fixed state variable. Prior results on the convex analytic methoc (in particular due to Borkar \cite{Borkar2} and nearly all the papers cited above \cite{kurano1989existence,hernandez1993existence,survey,HernandezLermaMCP,vega1999sample,lasserre1999sample,arapostathis2017some}) have assumed weak continuity of the kernel in both the state and action variables. Related to this contribution, recently \cite{yu2020minimum} established the existence of an optimal solution for countable action and Borel state spaces through majorization conditions via Lusin's theorem. A careful study of the topology of $w$-$s$ convergence,
which our existence analysis builds upon in this paper, reveals that Lusin's theory is
what establishes the connections between weak topology and the $w$-$s$ topology
via majorization conditions. Accordingly, in this paper the direct use of $w$-$s$ topology makes the analysis here more direct and concise, and as
opposed to the countable action space case (which makes functions continuous in the actions) in \cite{yu2020minimum}, here we consider general action spaces.
 
\item[(ii)] In Theorem \ref{optimalDeterministicStatHL}, we provide conditions under which the solution to
an optimal average cost stochastic control problem is a deterministic stationary policy. To our knowledge, there exists only two main such results
employing the convex analytic method, which as noted above generally require more
relaxed conditions compared with approaches directly utilizing the ACOE/ACOI. 
The first one is \cite[Proposition 9.2.5]{CTCN} and \cite[Lemma 2.4]{Borkar2}
for the countable probability space setup, and the second one due to \cite[Section 3.2]{Borkar2} for the continuous space setup, with the latter under restrictive conditions needed for applying Schauder's fixed point theorem. 
We also note that via a direct relationship between average cost optimality and ACOI and utilizing Blackwell \cite{Blackwell2,Blackwell3}; \cite[Corollary 5.4(b)]{hernandez1993existence} establishes the optimality of stationary and deterministic policies under a positive Harris recurrence assumption (see Section \ref{HLBound}), this analysis is utilized in Theorem \ref{optimalDeterministicStatHL}.

\item[(iii)] In some applications it may be useful to know not only that optimal policies
are deterministic, but that deterministic policies are dense in the sense of
approximability of the costs induced under randomized policies. In Theorem \ref{denseDet}, we establish conditions for not only the optimality, but also the
denseness of the attained performance values under deterministic (and also possibly quantized, i.e. with finitely many actions) stationary policies
in those attained under randomized stationary policies.
In other words, we show that, under mild conditions, the cost under any randomized
stationary policy can be approximated arbitrarily closely by the cost under
some deterministic stationary policy.
\end{itemize}

\section{The Convex Analytic Approach and a Refined Existence Result on the Optimality of Stationary (Possibly Randomized) Policies}

Recall that we are interested in the minimization
\begin{equation}\label{constOpt1}
\inf_{\gamma \in \Uadm} \limsup_{T \to \infty} \frac{1}{T} \Exp^{\gamma}_{x_0} \Biggl[\sum_{t=1}^T c(X_t,U_t)\Biggr]\,,
\end{equation}
where $\Exp^{\gamma}_{x_0}$ denotes the expectation over all sample paths with initial state given by $x_0$ under the admissible policy $\gamma$.

We refer the reader to \cite{ross1971nonexistence} for an example where an optimal policy may not be stationary under an average cost optimality criterion even for countable state/action spaces. Therefore, the conditions presented in the following are not superfluous.

\subsection{Some definitions}\label{Sdef}
We summarize here some definitions which we
use frequently in the paper.

For $\gamma\in\Usm$, we let
\begin{equation}\label{E-Tgamma}
\cT^\gamma (A\smid x) \,\df\, \int_{\cU(x)}
\cT(A\smid x,u)\,\gamma(\D{u}\smid  x)\,.
\end{equation}

We let $\cM_b(\XX)$ ($C_b(\XX)$)
 denote the space of bounded Borel measurable (continuous) real-valued functions
on $\XX$.
For $\mu\in\Pm(\KK)$ and $f\in\cM_b(\XX)$, we define
$\mu\cT\in\Pm(\XX)$ and $\cT f\colon \KK\to\RR$ by
\begin{equation}\label{E-muT}
\mu \cT(A) \,\df\, \int_{\KK} \mu(\D{x},\D{u}) \cT( A \smid x, u)\,,
\quad A\in\Bor(\XX)\,,
\end{equation}
and
\begin{equation}\label{E-cTf}
\cT f(x,u) \,\df\, \int_{\XX} f(y) \cT( \D{y} \smid x, u)\,,
\quad (x,u)\in\KK\,,
\end{equation}
respectively.

We use the convenient notation for integrals of functions
\begin{equation}\label{E-not1}
\mu(f)\,=\,\langle \mu, f\rangle \,\df\, \int_{\KK} f(x,u)\,\mu(\D{x},\D{u})\,,
\end{equation}
and similarly for $f\in\cM_b(\XX)$ and $\mu\in\Pm(\XX)$ if no ambiguity arises.
Clearly then, we have
\begin{equation*}
\langle \mu\cT, f\rangle\,=\, \langle \mu, \cT f\rangle
\qquad\text{for\ } \mu\in\Pm(\KK)\,,\ f\in\cM_b(\XX)\,.
\end{equation*}

The set of \emph{invariant occupation measures} (or, as is used more commonly in the literature: \emph{ergodic occupation measures}\footnote{It is perhaps more appropriate to use the term \emph{invariant occupation measures}, instead of \emph{ergodic occupation measures} since clearly the measures in $\cG$ are not all {\it ergodic}: we say that an invariant measure $\mu$ is ergodic if the support of $\mu$ does not contain two disjoint absorbing sets. However, traditionally the latter term has been used in the literature, see e.g. \cite{survey}.})
is defined by
\begin{equation*}
\cG \,\df\,
\bigl\{\mu\in\Pm(\KK)
\colon \mu(B \times\Act) = \mu \cT( B), \ B \in \Bor(\XX) \bigr\}\,.
\end{equation*}
We also let
\begin{equation*}
\cH \,\df\,
\biggl\{\uppi \in \Pm(\XX)
\colon \exists \gamma \in \Gamma_{S} \text{\ such that\ } \uppi(A) = \int_\XX
\cT^{\gamma}(A\smid x)\,\uppi(\D{x}),\  A \in \Bor(\XX) \biggr\}
\end{equation*}
denote the set of \emph{invariant probability measures} of the controlled
Markov chain.

Let $\mu\in\cG$.  It is well known that $\mu$ can be disintegrated into
a stochastic kernel $\phi$ on $\XX\times\Bor(\Act)$ and $\uppi\in\Pm(\XX)$ such that
\begin{equation*}
\mu(\D{x},\D{u}) \,=\, \phi(\D{u}\smid x)\,\uppi(\D{x})\,,
\end{equation*}
and $\phi$ is $\uppi$-a.e. uniquely defined on the support of $\uppi$.
We denote this disintegration by $\mu =\phi\circledast \uppi$.
Therefore, if $\gamma\in\Usm$ is any policy which agrees $\uppi$-a.e. with $\phi$,
then we have $\uppi(A) = \cT^{\gamma}(A\smid x)\,\uppi(\D{x})$
for $A\in\Bor(\XX)$.
Therefore, $\uppi\in\cH$.
Conversely, if $\uppi\in\cH$ with an associated $\gamma\in\Usm$, then it is clear from
the definitions that $\gamma\circledast\uppi\in\cG$.

Define
\begin{equation*}
\delta^* \,\df\, \inf_{\mu \in \cG}\, \langle \mu, c \rangle\,.
\end{equation*}
A measure $\mu\in\cG$ for which the infimum is attained is called \emph{optimal}.
\cref{S2.2,S2.3} concern the existence of optimal invariant occupation measures.

\subsection{Review: Optimality under weakly continuous kernels}\label{S2.2}

We first review the general proof method of some existing results, due to \cite{survey,Borkar2,kurano1989existence,hernandez1993existence,HernandezLermaMCP}, on the existence of an optimal $\mu\in\cG$ under the hypothesis that the transition kernel $\cT$ is weakly continuous.
This property is defined as follows.

\begin{itemize}
\item[\hypertarget{H1}{\textbf{(H1)}}]
The transition kernel $\cT$ is called \emph{weakly continuous} if the map
\begin{equation*}
\KK\,\ni\,(x,u)\mapsto\int_{\XX} f(z)\cT(\D{z}\smid x,u)
\end{equation*}
is continuous for all $f\in C_b(\XX)$.
\end{itemize}

Continuing, for $T \ge 1$, we let
\begin{equation*}
v_T(D) \,=\, \frac{1}{T} \sum_{t=0}^{T-1} \Ind_D(X_t,U_t),
\quad D \in \Bor(\XX\times\Act)\,.
\end{equation*}
Consider any policy $\gamma$ in $\Uadm$, $X_0 \sim \nu$, and let for $T \ge 1$,
\begin{equation*}
\mu_T^\gamma(D)  \,=\, \Exp_{\nu}^{\gamma}[v_T(D)] \,=\,  \frac{1}{T} \Exp^{\gamma}_{\nu} 
\Biggl[ \sum_{t=0}^{T-1} \Ind_D(X_t,U_t) \Biggr]\,,
\quad D \in \Bor(\XX\times\Act)\,.
\end{equation*}
We refer to $\bigl\{\mu_T^\gamma\bigr\}_{T>0}$ as the family of
\emph{mean empirical occupation measures}
under the policy $\gamma\in\Uadm$, and with initial distribution $\nu$.
Through what is often referred to as a \emph{Krylov-Bogoliubov-type} argument,
for every $A\in\Bor(\XX)$, we have
\begin{equation}\label{KBarg}
\begin{aligned}
\babs{\mu_{T}^\gamma(A\times\Act) - \mu_{T}^\gamma \cT (A)}
& \,=\, \frac{1}{T}\,
\Babss{ \Exp^{\gamma}_{\nu}\Biggl[ \sum_{t=0}^{T-1} \Ind_{A\times\Act}(X_t,U_t)
- \sum_{t=1}^{T} \Ind_{A\times\Act}(X_t,U_t) \Biggr]}\\
& \,\le\, \frac{1}{T} \,\to\, 0\quad\text{as\ } T\to\infty\,.
 \end{aligned}
 \end{equation} 
Observe that \eqref{KBarg} holds for any policy $\gamma \in \Uadm$.

Suppose that, along some subsequence $\{t_k\}\subset\NN$,
$\mu^\gamma_t$ converges weakly to some $\mu\in\Pm(\KK)$, which we denote
as $\mu^\gamma_{t_k}\Rightarrow \mu$.
Using \cref{E-not1}, we
write the triangle inequality
\begin{equation}\label{E-triangl}
\begin{aligned}
\babs{\mu(f) - \mu\cT (f)}
&\,\le\,  \babs{\mu(f) - \mu^\gamma_{t_k}(f)}
 + \abs{\mu^\gamma_{t_k}(f) - \mu^\gamma_{t_k}\cT(f)} \\[3pt]
 &\mspace{100mu} +\babs{\mu^\gamma_{t_k}\cT(f) - \mu\cT (f)}
\end{aligned}
\end{equation}
for $f\in C_b(\XX)$.
This notation is consistent since $f$ may be viewed also as an element of $C_b(\KK)$.
Suppose that \hyperlink{H1}{(H1)} holds.
The first term on the right hand side of \cref{E-triangl} vanishes as $k\to\infty$ by weak convergence,
while the second term does the same by \cref{KBarg}.
Since 
\begin{eqnarray}\label{weakFC}
\mu^\gamma_{t_k}\cT(f)=\mu^\gamma_{t_k}(\cT f),
\end{eqnarray}
 and $\cT f\in C_b(\KK)$ by
\hyperlink{H1}{(H1)}, it follows that the third term also vanishes
as $k\to\infty$ by the weak convergence $\mu^\gamma_{t_k}\Rightarrow \mu$.
Since the class $C_b(\XX)$ distinguishes points in $\Pm(\XX)$,
this shows that $\mu(A,\Act)=\mu\cT(A)$ for all $A\in\Bor(\XX)$, which implies that
$\mu\in\cG$ by the definition of the latter.
Thus we have shown the following.

\begin{lemma}\label{L2.1}
Under \hyperlink{H1}{(H1)}, the limit
of any weakly converging subsequence of mean empirical occupation measures
is in $\cG$.
\end{lemma}

Recall (\ref{E-AC}). This expected cost can be equivalently written as
\begin{equation*}
J(x,\gamma)
\,\df\, \limsup_{T \to \infty}\, \langle \mu_T^\gamma, c \rangle\,,
\end{equation*}
where $\mu_T^\gamma$ is the mean empirical occupation measure under $\gamma$.
Let $\{t_k\}\subset\NN$ be a subsequence along which
$\langle \mu_{t_k}^\gamma, c \rangle$ converges to $J(x,\gamma)$
and suppose that $\mu_{t_k}\Rightarrow \mu\in\cG$.
Then
\begin{equation}\label{E-thmWeak0}
J(x,\gamma)\,=\,\liminf_{t_k \to \infty}\, \langle \mu_{t_k}^\gamma, c \rangle 
\,\ge\, \Bigl\langle \lim_{t_k \to \infty}\,  \mu_{t_k}^\gamma, c \Bigr\rangle
\,=\, \langle\mu, c\rangle
\,\ge\,
\delta^*\,,
\end{equation}
where for the first inequality we use the fact
that, since $c$ is lower semi-continuous (l.s.c.)
and bounded from below, the map $\mu\to \langle \mu,c\rangle$ is lower semi-continuous.
The above shows that $J^*(x) \,\ge\, \delta^*$. We now establish conditions for which the above is indeed an equality.


\begin{assumption}\label{A2.2}
\begin{itemize}
\item[(A)] The state  and action spaces $\XX$ and $\Act$ are Polish.
The set-valued map $\cU\colon\XX\to\Bor(\Act)$ is upper semi-continuous
and closed-valued.
\item[(A')] The state  and action spaces $\XX$ and $\Act$ are compact.
The set-valued map $\cU\colon\XX\to\Bor(\Act)$ is upper semi-continuous
and closed-valued.
\item[(B)]
The non-negative
running cost function $c(x,u)$ is  is l.s.c. and $c\colon\KK\to\RR$ is inf-compact, i.e. $\{(x,u) \in \KK: c(x,u) \leq \alpha\}$ is compact for every $\alpha \in \mathbb{R}_+$.
\item[(B')] The cost function $c$ is bounded and l.s.c..
\item[(C)]
There exists a policy and an initial state leading to a finite cost $\eta \in \mathbb{R}_+$.
\item[(D)] \hyperlink{H1}{(H1)} holds.
\item[(E)] Under every stationary policy, the induced Markov chain is Harris recurrent.
\end{itemize}

\end{assumption}

Before we present a theorem, we now review the following concerning ergodic properties of (control-free) Markov chains: Let $c \in L_1(\mu):= \{f: \mathbb{X} \to \mathbb{R}, \int |f(x)| \mu(dx) < \infty \}$. Suppose that $\mu$ is an invariant ergodic probability measure for an $\mathbb{X}$-valued Markov chain. Then, it follows that for $\mu$ almost everywhere $x \in \mathbb{X}$:
\[\lim_{T \to \infty} {1 \over T}\sum_{t=1}^T c(X_t) = \int c(x) \mu(dx),\]
$P_x$ almost surely (that is conditioned on $x_0=x$, with probability one, the above holds). Furthermore, again with $c \in L_1(\mu)$, for $\mu$ almost everywhere $x \in \mathbb{X}$
\begin{eqnarray}\label{IETMean}
\lim_{T \to \infty} {1 \over T} E_x\bigg[ \sum_{t=1}^T c(X_t) \bigg] = \int c(x) \mu(dx),
\end{eqnarray}
On the other hand, the positive Harris recurrence property allows the almost sure convergence to take place for every initial condition: If $\mu$ is the invariant probability measure for a {\it positive Harris recurrent} Markov chain, it follows that for all $x \in \mathbb{X}$ and for every $c \in L_1(\mu)$ \cite[Theorem 17.1.7]{MeynBook} or \cite[Theorem 4.2.13]{HernandezLasserreErgodic}
\begin{eqnarray}
\lim_{T \to \infty} {1 \over T}\sum_{t=1}^T c(X_t) = \int c(x) \mu(dx),\label{ergodicPHR1}
\end{eqnarray}
almost surely.
However, for every $c \in L_1(\mu)$, while (\ref{ergodicPHR1}) holds for all $x \in \mathbb{X}$, it is not generally true that \cite[Chapter 14]{MeynBook} (see e.g. \cite[Example 3.1]{yu2020minimum}) that 
\[\lim_{T \to \infty} {1 \over T}E_x[\sum_{t=1}^T c(X_t)] = \int c(x) \mu(dx), \]
for all $x \in \mathbb{X}$. Thus, we can not in general relax the boundedness condition for the convergence of the expected costs. 
However, with $c$ bounded, for all $x \in \mathbb{X}$
\begin{eqnarray}\label{convPHRBounded}
\lim_{T \to \infty} {1 \over T} E_x\bigg[ \sum_{t=1}^T c(X_t) \bigg] = \int c(x) \mu(dx)
\end{eqnarray}
This follows as a consequence of Fatou's lemma and (\ref{ergodicPHR1}). Further refinements are possible via return properties to small sets and $f$-regularity of cost functions \cite{ABor-19, MeynBook} (e.g. this convergence holds if \cite[Theorem 14.0.1]{MeynBook} holds and $X_0=x$ with $x \in \{z: V(z) < \infty\}$). We refer the reader to \cite[Chapters 14 and 17]{MeynBook} or \cite[Chapters 2 and 4]{HernandezLasserreErgodic} for additional discussions.


\begin{theorem}\label{thmWeak}
a) Under \cref{A2.2}\,\ttup{A, B, C, D} there exists an optimal measure in $\cG$.
b) Under \cref{A2.2}\,\ttup{A', B', D, E}, there exists a policy in $\Usm$ which is optimal for the control problem given in \cref{constOpt1} for every initial condition. 
\end{theorem}

\begin{proof}

a) 
Consider \cref{A2.2}\,\ttup{A, B, C, D}. By (B, C) we have that the set of policies $\gamma$ which lead to a bounded cost is so that $\langle \mu^\gamma_T,c \rangle \leq M < \infty$ for all $T$, which implies that $\{\mu_T^\gamma\,,\, T>0\}$ is tight. Thus along some subsequence $\mu_{t_k}\to\mu\in\Pm(\KK)$. 
As shown in the paragraph preceding \cref{L2.1}, $\mu\in\cG$. 

Furthermore, under hypothesis (A), the set $\mathbb{K}= \cup_{x \in \XX} \{(x,u), u \in \cU(x)\}$ is closed by \cite[Lemma D.3]{HernandezLermaMCP}. Thus, by the Portmanteau theorem every weak limit of a converging sequence of probability measures on $\mathbb{K}$ is also supported on $\mathbb{K}$.

Consider a sequence $\{\mu_k\}_{k\in\NN}\subset\cG$ such that
$\langle\mu_k,c\rangle \to \delta^*$ as $k\to\infty$, the sequence $\mu_{t_k}$ is tight by inf-compactness, and any limit point $\mu_*$ of this sequence
is in $\cG$ with $\mu_*(\mathbb{K})=1$. Thus, by \cite[Prop. D.8]{HernandezLermaMCP} we have an optimal control policy $\phi$. Taking limits as in \cref{E-thmWeak0}, we obtain $\langle\mu_*,c\rangle=\delta^*$. This establishes the first part of the theorem.

Define a stationary policy $\gamma$ via the disintegration
\begin{equation}\label{denklem11defineP}
 \mu_*(\D{x},\D{u})\,=\,\gamma_*(\D{u} \smid x)\,\uppi_*(\D{x})
 \end{equation}
$\mu_*$ almost surely.
Note that via this disintegration the control $\gamma_*$ is defined
$\uppi_*$-a.e. Take $\phi\in\Usm$ be any policy that agrees with $\gamma_*$ on the support of $\uppi_*$. 

b) Under (A', B', D), via (\ref{weakFC}) and that $\cT f\in C_b(\KK)$ by \hyperlink{H1}{(H1)}, we have that $\cG$ is compact; we also have that Portmanteau theorem applies as in part a). By hypothesis (E), since the chain under an optimal $\phi$, is Harris recurrent, $\uppi_*$ is its unique invariant probability measure. Optimality of $\phi$, for every initial condition, then follows by positive Harris recurrence given that $c$ is bounded via (\ref{convPHRBounded}) under hypothesis (B').
Thus, $J(x,\phi) = \langle\mu_*,c\rangle$ and $\mu_*$ is optimal. 
\end{proof}

\cref{thmWeak} can be stated under weaker assumptions.
See, for example, \cite[Theorem~2.1]{ABor-19} among other references in the literature.
We have chosen to state it under somewhat stronger hypotheses
in order to present a simple and short proof that conveys
the essential arguments. 

In general, in the absence of \cref{A2.2}\,(E),
there is a consideration of reachability.
Suppose that the chain under the policy $\phi$ as defined in the proof
of \cref{thmWeak} is a $T$ model (see \cite{Tweedie-94}).
Then, as asserted in \cite[Theorem~6.1]{Tweedie-94}, the Doeblin decomposition of the state space
contains, in general, a countable collection of maximal Harris sets.
In particular, we have a decomposition
into the disjoint union $\XX \,=\, \bigl(\cup_{i\in\NN} H_i\bigr)\cup E\,,$ where each $H_i$ is a maximal Harris set with invariant measure $\uppi_i$, and $E$ is transient.
Now, by part (ii) of Theorem~6.1 in \cite{Tweedie-94}, only a finite number of
the sets $H_i$ may have a nonempty intersection with any given compact set. Thus, if the Markov Chain is not
recurrent, the stationary policy defined above, in general,
is only optimal in a restricted set of initial conditions. On implications related to insensitivity to such initial state dependence, the reader is referred to \cite{lasserre1999sample} and \cite[Prop. 11.4.4(c) and Lemma 11.4.5(a)]{hernandezlasserre1999further}, among other references, for further results on sample path average cost optimality and expected average cost optimality.


\subsection{New conditions: optimality under setwise convergence and strong continuity in actions for each state}\label{S2.3}

There are many important applications where the kernel $\cT$ is not weakly continuous.
For example, consider dynamics described by a stochastic difference equation
on $\RR^d$ of the form
\begin{equation*}
X_{n+1}\,=\, F(X_{n},U_{n})+W_{n}\,,\qquad n=0,1,2,\dotsc\,,
\end{equation*}
where $\XX=\mathbb{R}^{n}$ and the $W_{n}$'s are independent and
identically distributed (i.i.d.) random vectors whose distribution
has a bounded and continuous density function. We assume that $F$ is bounded and $u\mapsto F(x,u)$ is continuous
for all $x\in\XX$. It is clear that the transition kernel $\cT$ is not, in general, weakly
continuous. However, it satisfies the following hypothesis.

\smallskip
\begin{itemize}
\item[\hypertarget{H2}{\textbf{(H2)}}]
The transition kernel $\cT$ satisfies the following:
\begin{itemize}
\item[(a)]
For any $x\in\XX$, the map $u\mapsto\int f(z)\cT(\D{z}\smid x,u)$ is continuous
for every bounded measurable function $f$.
\item[(b)] 
There exists a finite measure $\nu$ majorizing $\cT$, that is
\begin{equation}\label{EH2A}
\cT(\D{y} \smid x, u) \,\le\, \nu(\D{y})\,, \qquad x \in \XX,\; u \in \Act\,.
\end{equation}
\end{itemize}
\end{itemize}




\begin{assumption}\label{A2.5}
The following hold:
\begin{enumerate}
\item[(A)] The state and action spaces $\XX$ and $\Act$ are Polish. The set $\mathbb{K}= \cup_{x \in \XX} \{(x,u), u \in \cU(x)\}$ is measurable (see \cite[Lemma D.3]{HernandezLermaMCP} for conditions) and the set-valued map $\cU\colon\XX\to\Bor(\Act)$ is compact-valued.
\item[(A')] The state  and action spaces $\XX$ and $\Act$ are compact.
The set $\mathbb{K}$ is measurable and set-valued map $\cU\colon\XX\to\Bor(\Act)$ is compact-valued.
\item[(B)]
The non-negative
running cost function $c(x,u)$ is continuous in $u\in\cU(x)$ for every $x\in\XX$
and $c\colon\KK\to\RR$ is inf-compact.
\item[(B')] The cost function $c$ is bounded, and continuous in $u\in\cU(x)$ for every $x\in\XX$.
\item[(C)]
There exists a policy and an initial state leading to a finite cost $\eta \in \mathbb{R}_+$.
\item[(D)] \hyperlink{H2}{(H2)} holds.
\item[(E)] Under every stationary policy, the induced Markov chain is Harris recurrent.
\end{enumerate}
\end{assumption}

Let us recall the $w$-$s$ topology studied by Sch\"al \cite{schal1975dynamic} (see Balder \cite{balder2001} for further properties).

\begin{definition}
The $w$-$s$ topology on the set of probability measures $\mathcal{P}(\XX \times \Act)$
is the coarsest topology under which
$\int f\,\D\nu\colon \mathcal{P}(\XX \times \Act) \to \mathbb{R}$
is continuous for every measurable and bounded $f(x,u)$ which is continuous in $u$
for every $x$
(but unlike the weak topology, $f$ does not need to be continuous in $x$).
\end{definition}

It is a consequence of \cite[Theorem 3.10]{schal1975dynamic} or
\cite[Theorem 2.5]{balder2001} that \cref{EH2A}, by implying setwise sequential
pre-compactness of marginal measures on the state, ensures that every weakly
converging sequence of mean
empirical occupation measures also converges in the $w$-$s$ sense.
\Cref{EH2A} implies setwise sequential pre-compactness by
\cite[Proposition 3.2]{SaldiLinderYukselTAC14}, which in turn builds on
\cite[Corollary 1.4.5]{HernandezLasserreErgodic}; see also \cite[Theorem 4.17]{Hai06}.


\begin{theorem}\label{thmSetwise}
a) Under \cref{A2.5}\,\ttup{A, B, C, D}, there exists an optimal measure in $\cG$.
b) Under \cref{A2.5}\,\ttup{A', B', D, E}, there exists a policy in $\Usm$ which is optimal for the control problem given in \cref{constOpt1} for every initial condition.
\end{theorem}

First, note the following counterpart to Lemma \ref{L2.1}.
\begin{lemma}\label{L2.2}
Under \hyperlink{H2}{(H2)}, the limit
of any $w$-$s$ converging subsequence of mean empirical occupation measures
is in $\cG$.
\end{lemma}

\begin{proof}
We follow the notation used in the discussion leading to Lemma \ref{L2.1}. Suppose that, along some subsequence $\{t_k\}\subset\NN$,
$\mu^\gamma_t$ converges to some $\mu\in\Pm(\KK)$ in the $w$-$s$ sense, which we denote
as $\mu^\gamma_{t_k}\Rightarrow \mu$.
As in (\ref{E-triangl}) we have the triangle inequality
\begin{equation}\label{E-triangl2}
\begin{aligned}
\babs{\mu(f) - \mu\cT (f)}
&\,\le\,  \babs{\mu(f) - \mu^\gamma_{t_k}(f)}
 + \abs{\mu^\gamma_{t_k}(f) - \mu^\gamma_{t_k}\cT(f)} \\[3pt]
 &\mspace{100mu} +\babs{\mu^\gamma_{t_k}\cT(f) - \mu\cT (f)}
\end{aligned}
\end{equation}
for $f\in {\cal M}_b(\XX)$.
If \hyperlink{H2}{(H2)} holds, the first term on the right hand side of \cref{E-triangl2} vanishes as $k\to\infty$ by $w$-$s$ convergence,
while the second term does so by \cref{KBarg}.
We have 
\begin{eqnarray}\label{wsFC}
\mu^\gamma_{t_k}\cT(f)=\mu^\gamma_{t_k}(\cT f),
\end{eqnarray}
where $\cT f$ is as defined in \cref{E-cTf}. Since  $\cT f$ is continuous in $u$ for every fixed $x$, by \hyperlink{H2}{(H2)}, it follows that the third term also vanishes
as $k\to\infty$ by the $w$-$s$ convergence $\mu^\gamma_{t_k}\Rightarrow \mu$.
This shows that $\mu(A,\Act)=\mu\cT(A)$ for all $A\in\Bor(\XX)$, which implies that
$\mu\in\cG$.
\end{proof}

\begin{proof}[Proof of \cref{thmSetwise}]

The proof follows along the lines of \cref{thmWeak},
but instead of weak convergence, we work with $w$-$s$ convergence.

As noted earlier, \cref{EH2A} ensures that every weakly converging sequence also does
so under the $w$-$s$ sense (see \cite[Proposition 3.2]{SaldiLinderYukselTAC14},
which in turn builds on \cite[Corollary 1.4.5]{HernandezLasserreErgodic} or
\cite[Theorem 4.17]{Hai06}). 

a) Accordingly,
\cref{A2.5}\,\ttup{A, B, C, D} ensures that
each mean empirical occupation measure leading to a finite cost has a weakly converging subsequence, and which then is a $w$-$s$ converging subsequence. Lemma \ref{L2.2} then implies that the limit of this sequence $\mu$ is in ${\cal G}$. 

Furthermore, under hypothesis (A) or (A'), the set $\mathbb{K}$ is measurable. Thus, by the generalized Portmanteau theorem \cite[Proposition 3.2]{balder2001} for $w$-$s$ convergence every $w$-$s$ limit of a converging sequence of probability measures $\mu_n$ with $\mu(\mathbb{K})=1$ is also supported on $\mathbb{K}$.


Now, \Cref{EH2A} implies also that the set of measures in $\cG$ leading to a cost less than $\eta$ is $w$-$s$
pre-compact, that is, for every sequence $\mu_n \in \cG$ with $\langle \mu_n, c \rangle \leq \eta$,
 there exists a subsequence which converges (in the $w$-$s$ sense) to a limit: Now, let $\mu_n$ be a sequence in $\cG$
such that $\mu_n\xrightarrow[]{w\text{-}s}\mu$. We show that $\mu \in \cG$ and this also leads to a cost less than $\eta$. 

Using the definition in \cref{E-muT}, we note first that
\begin{equation*}
\begin{aligned}
\cG &\,=\, \biggl\{\mu \in \mathcal{P}(\XX\times\Act)\,
\colon  \int_{\XX\times\Act} f(x) \mu(\D{x},\D{u})
\,=\,  \int_{\XX} f(y)\,\mu\cT(\D{y})\,, \quad \forall\,f \in \cM_b(\XX) \biggr\},
\end{aligned}
\end{equation*}
where as defined in \cref{Sdef}, $\cM_b(\XX)$ denotes the set of bounded
Borel measurable functions on $\XX$. 
Thus, for every $f \in \cM_b(\XX)$, we have
\begin{equation}\label{PT2.6A}
\begin{aligned}
\lim_{n \to \infty} \int_{\XX\times\Act} f(x) \mu_n(\D{x},\D{u})
&\,=\, \lim_{n \to \infty} \int_{\XX} f(y)\,\mu_n\cT(\D{y})\\
&\,=\, \lim_{n \to \infty} \int_{\XX} \cT f(x,u) \mu_n(\D{x},\D{u})\,,
\end{aligned}
\end{equation}
where $\cT f$ is as defined in \cref{E-cTf}.
Since  $\cT f$ is continuous in $u$ for every fixed $x$,
by \cref{A2.5}\,(D), and $\mu_n\xrightarrow[]{w\text{-}s}\mu$, we obtain
\begin{equation*}
\lim_{n \to \infty}\, \langle \mu_n,\cT f\rangle
\,=\,  \langle \mu,\cT f\rangle
\,=\, \langle \mu\cT, f\rangle\,,
\end{equation*}
and
\begin{equation*}
\lim_{n \to \infty}\, \langle \mu_n, f\rangle
\,=\,  \langle \mu, f\rangle\,.
\end{equation*}
Since the terms on the left hand side are equal  by \cref{PT2.6A},
we have equality of the terms on the right hand side, which
implies that $\mu \in \cG$.

Note that the integral $\langle \mu,c\rangle$ is lower semi-continuous in $\mu$.
This follows
by truncating $c$ as $c^N(x,u) = \min(N, c(x,u))$, and then taking the limit $N \to \infty$ noting that for every finite $N$, $\langle\mu, c^N\rangle$ is continuous in $\mu$ by
the $w$-$s$ convergence. Thus, we also have that $\langle \mu,c\rangle \leq \eta$.
As a result, there exists an optimal measure $\mu_* \in \cG$ with $\mu_*(\mathbb{K})=1$, and by, e.g., \cite[Prop. D.8]{HernandezLermaMCP}, we have an optimal control policy $\phi$. 

b) Now, under (A') $w$-$s$ compactness follows from the existence of a $w$-$s$ converging subsequence and (\ref{wsFC}) and the discussion following it. Then, under \cref{A2.5}\,\ttup{A', B', D} and \ttup{E}, as in \cref{thmWeak}, optimality of $\phi$ for every initial condition follows by positive Harris recurrence given that $c$ is bounded via (\ref{convPHRBounded}).
\end{proof}

%

\section{Optimality of Deterministic Stationary Policies}

In this section, we  provide conditions under which an optimal average cost stochastic control problem is a \emph{deterministic} policy.

\subsection{Preliminaries}

\begin{definition}\label{defDetRandPol}
A policy $\gamma\in\Usm$ under which the chain has an invariant probability measure
$\uppi_\gamma$, is called $\uppi_\gamma$-deterministic (or simply, deterministic), if
\begin{equation*}
\uppi_\gamma\bigl(\{x\in\XX\colon \gamma(\,\cdot\smid x) \mathrm{\ is\ Dirac}\}\bigr)
\,=\,1\,.
\end{equation*}
If the policy is not $\uppi_\gamma$-deterministic, we say that
it is $\uppi_\gamma$-randomized (or simply, randomized).
\end{definition}

Here, $\mu_\phi$ denotes the invariant occupation measure of
the chain under a stationary Markov policy $\phi$.

%
%

\subsubsection{Convexity of the set of invariant occupation measures}
Under the conditions presented in the previous section, the space $\cG$ is closed under either the weak convergence or the $w$-$s$ topologies. 

We now discuss convexity of $\cG$. Let $\kappa \in (0,1)$ and consider two
invariant occupation measures $\mu^1, \mu^2 \in \cG$. Let
\begin{equation}\label{E3.1A}
\mu^i(\D{x},\D{u})= \phi^i(\D{u}\smid x) \uppi_{\phi^i}(\D{x})\quad\text{for\ }i=1,2\,,
\end{equation}
denote their disintegration into invariant probability measures $\uppi_{\phi^i}$,
and Markov policies $\phi^i$, $i=1,2$, respectively.
Define
\begin{equation}\label{E3.1B}
\uppi(\D{x}) \,\df\, \kappa \uppi_{\phi^1}(\D{x}) + (1 - \kappa)\uppi_{\phi^2}(\D{x})\,.
\end{equation}
Note that $\uppi(\D{x}) = 0 \implies \uppi_{\phi^i}(\D{x})=0$ for $i=1,2$.
As a consequence, the Radon-Nikodym derivative of $\uppi_{\phi^i}$ with respect to $\uppi$ exists.
Let $f^i(x)\df\frac{\D\uppi_{\phi^i}}{\D\uppi}(x)$, $i=1,2$, and
\begin{equation}\label{E3.1C}
\phi(\D{u}\smid x) \,\df\, \kappa f^1(x) \,\phi^1(\D{u}\smid x)
+ (1-\kappa) f^2(x)\,\phi^2(\D{u}\smid x)\quad \uppi\text{-a.e.}
\end{equation}
Then
\begin{equation}\label{E3.1D}
\mu(\D{x},\D{u})\,\df\,
\phi(\D{u}\smid x)\uppi(\D{x}) = \kappa \mu^1(\D{x},\D{u}) + (1-\kappa) \mu^2(\D{x},\D{u})\,,
\end{equation}
and it follows by applying the definition
that $\mu\in\cG$.
Therefore, $\cG$ is convex.

In the following we let $\cG_e$ denote the set of extreme points of $\cG$.

\subsubsection{A partial characterization of \texorpdfstring{$\cG_e$}{}}

\begin{lemma}
If a measure $\mu$ is not in $\cG_e$, then one of the following conditions are satisfied:
(i) The control policy inducing it is randomized, or
(ii) under this policy the Markov chain has multiple invariant probability measures.
\end{lemma}

\begin{proof} 
Let $\mu$ be an invariant occupation measure in $\cG$ which is not extreme.
This means that there exist $\kappa \in (0,1)$ and distinct
invariant occupation measures $\mu^1, \mu^2 \in\cG$ such that
\cref{E3.1A,E3.1B,E3.1C,E3.1D} hold.

Suppose $f^1f^2=0$ $\uppi$-a.e.
Then the invariant measures $\uppi_{\phi^i}$, $i=1,2$,
 are singular with respect to each other, so under the policy $\phi$ the Markov chain has two distinct invariant probability measures. On the other hand, if $f^1f^2 \neq 0$ on a set of positive $\uppi_{\phi^i}$ measure, then by \cref{E3.1C}, the policy is randomized on that set.
\end{proof}


However, the converse direction is more consequential for optimization purposes,
as we wish to show the optimality of deterministic policies.
Towards this end, in what follows, we
characterize the extreme points of the convex set $\cG$.
Since an optimal solution can,
without any loss of generality, be searched over the extreme points of this set
due to the linear programming formulation, this characterization provides
insights on the structure of optimal policies.
In particular, we establish the optimality of deterministic stationary policies.

\subsubsection{Revisiting the countable state/action space setup: Optimality of deterministic policies}

As noted earlier, the countable setup has been studied in \cite[2.4]{Borkar2} and \cite[Proposition 9.2.5]{CTCN}. We provide a different proof for \cref{extremePoints} which may also be utilized in the continuous space setup, see \cref{S3.5}. 



Following Definition \ref{defDetRandPol}, if $\phi$ is a non-deterministic policy, we can select $\alpha\in\XX$ and lying on
the support of $\mu_\phi$, such that $\phi(\D{u}\smid \alpha)$ can be
expressed as a non-trivial convex combination of two different probability measures $ \gamma_1$ and $\gamma_2$ on $\Act$
\begin{equation}\label{randomizedP}
\phi(\D{u}\smid \alpha) \,=\,\theta \gamma_1(\D{u}) + (1-\theta)\gamma_2(\D{u})\,,
\end{equation}
and $\theta\in(0,1)$.

\begin{lemma}\label{extremePoints}
We assume that the chain is controlled by some $\phi\in\Usm$
has an invariant probability measure $\uppi_\phi$.
Suppose that $\phi$ is non-deterministic on some set
 that has positive $\uppi_\phi$ measure.
Then the corresponding invariant occupation measure $\mu_\phi$
cannot lie in $\cG_e$.
\end{lemma}

\begin{proof}
Let  $\phi$ be a non-deterministic policy so that \cref{randomizedP} holds.
Let $\phi^i$, $i=1,2$, denote the Markov policy which at $\alpha$ (with $\uppi_\phi(\alpha) >0$) selects an action under $\gamma_i$ and agrees with $\phi$ everywhere else.
It is clear that, with $\tau_{\alpha}=\min(k>0: x_k=\alpha)$ denoting
the first return time to $\alpha$, we have the stochastic representations
\begin{equation}\label{L3.1A}
\uppi_{\phi^i}(x) \,=\, \frac{\Exp_{\alpha}^{\phi^1} \Bigl[\sum_{k=0}^{\tau_{\alpha}-1}
\Ind_{\{X_k=x\}}\Bigr]}{\Exp_{\alpha}^{\phi^1}[\tau_{\alpha}]}\,,\quad
i=1,2\,,
\end{equation}
and
\begin{equation*}
\begin{aligned}
\uppi_\phi(x) &\,=\,
\frac{\theta \Exp_\alpha^{\phi^1}
\Bigl[\sum_{k=0}^{\tau_{\alpha}-1} \Ind_{\{X_k=x\}}\Bigr] 
+ (1-\theta) \Exp_\alpha^{\phi^2}
\Bigl[\sum_{k=0}^{\tau_\alpha-1} \Ind_{\{X_k=x\}}\Bigr]}
{\theta\Exp_{\alpha}^{\phi^1}[\tau_\alpha] 
+  (1-\theta)\Exp_\alpha^{\phi^2}[\tau_{\alpha}]}\\
&\,=\, \kappa \uppi_{\phi^1} + (1-\kappa)\uppi_{\phi^2}\,,
\end{aligned}
\end{equation*}
where in the second equality we use \cref{L3.1A}, and the constant
$\kappa\in(0,1)$ defined by
\begin{equation*}
\kappa\,\df\, \frac{\theta\Exp_{\alpha}^{\phi^1}[\tau_\alpha] }
{\theta\Exp_{\alpha}^{\phi^1}[\tau_\alpha] 
+  (1-\theta)\Exp_\alpha^{\phi^2}[\tau_{\alpha}]}\,.
\end{equation*}
It follows from \cref{E3.1D} that
$$\mu_\phi = \kappa \phi^1\circledast \uppi_{\phi^1}
+ (1-\kappa)\phi^2\circledast \uppi_{\phi^2}\,.$$
It is clear that $\uppi_{\phi^i}(\alpha)>0$ for $i=1,2$.
Thus $\phi^1\circledast \uppi_{\phi^1}\ne\phi^2\circledast \uppi_{\phi^2}$
since the $\gamma_i$'s are not identical.
This shows that $\mu_\phi\notin\cG_e$.
\end{proof}

As a result, we can deduce that for such countable state and action spaces an optimal policy is stationary \emph{and deterministic}, provided that the convex analytic method is applicable.

\subsection{Uncountable standard Borel setup: Optimality of deterministic policies}

For an uncountable setup, the optimality of deterministic policies under the convex analytic approach has been an open problem with partial results available. In the following, we both present a review of relevant results and present further conditions.

\subsubsection{Arriving at ACOE/ACOI from the convex analytic method}\label{HLBound}

In general, establishing conditions for the existence of a solution to ACOI is an unfinished problem. Our findings reported earlier through the convex analytic method, via a duality analysis, may provide further conditions. One may express the linear program
\begin{equation}\label{PLPInf}
\min_{\nu\in\cG}\,\langle \nu, c \rangle
\end{equation}
as an infinite dimensional linear program, present its convex dual formulation and arrive at the ACOI. This then would lead to an existence result. 

However, a more direct argument (without using duality) along this direction was presented by Hern\'andez-Lerma \cite[Theorem 5.3]{hernandez1993existence}. This result shows that an average cost optimal randomized policy $\phi$, with invariant measure $\uppi_{\phi}$ satisfies the ACOI $\pi_{\phi}$ almost everywhere:
\begin{align}\label{ACOI11}
g + h(x) \geq  c(x,\phi(x)) + \int h(x') {\cal T}(dx'|x, \phi(x)) 
\end{align}
where $h$ is bounded from below. If one can ensure that the above holds for all $x \in \mathbb{X}$ (and not just $\pi_{\phi}$ almost everywhere) \cite[Prop. 5.2]{hernandez1993existence} shows that under this condition on $h$, (\ref{ACOI11}) implies that such a policy is indeed optimal. Again, if the above holds for all $x \in \mathbb{X}$, by utilizing Blackwell's theorem of optimality of deterministic policies (also called irrelevant information theorem) \cite{Blackwell2,Blackwell3}, we can replace $\phi$ with a deterministic $f \in \Gamma_{SD}$, which will then be optimal \cite[Corollary 5.4(b)]{hernandez1993existence}.

In the following, we relax the condition of (\ref{ACOI11}) holding for every $x$. The following is a refinement on \cite[Prop. 5.2]{hernandez1993existence}.

Let $g$ be a constant and $h : \mathbb{X} \to \mathbb{R}_+, f: \mathbb{X} \to {\cal P}(\mathbb{U})$ be so that for all $x \in B$ for some Borel set $B \subset \mathbb{X}$, 
\begin{align}\label{ACOI1a}
g + h(x) \geq \bigg(c(x,f(x) + \int h(x') {\cal T}(dx'|x,f(x)) \bigg) := \int \bigg(c(x,u) + \int h(x') {\cal T}(dx'|x,u) \bigg) f(du|x)
\end{align}

\begin{lemma}\label{ACOIEqn}
Let (\ref{ACOI1a}) hold with
\begin{eqnarray}
\liminf_{n \to \infty} {1 \over n}E^{\gamma^*}_x[h(X_n)] \geq 0, \label{condConv0}
\end{eqnarray}
for all $x \in B$ where $\gamma^* = \{f,f,f,\cdots\}$ and $P^{\gamma^*}(x,B) = 1$ for all $x \in B$. Then the stationary (possibly randomized) policy $\gamma^* = \{f,f,f,\cdots\}$ satisfies
\[g \geq J(x,\gamma^*),\]
for all $x \in B$.
\end{lemma}

\begin{proof}
We have
\begin{eqnarray}
&& E^{\gamma}[h(X_t) | x_{[0,t-1]},u_{[0,t-1]}] = \int_y h(y) P(X_t \in dy|x_{t-1},u_{t-1}) \\
&& \quad \quad \quad = c(x_{t-1},u_{t-1}) + \int_y h(y) P(dy|x_{t-1},u_{t-1}) - c(x_{t-1},u_{t-1}) 
\end{eqnarray}
By iterated expectations,
\[E^{\gamma^*}_x\bigg[\sum_{t=1}^n  h(X_t) - E^{\gamma^*}[h(X_t) | X_{[0,t-1]},U_{[0,t-1]}] \bigg] =0 \]
Now, under $\gamma^*$ we have that $B$ is an absorbing set and thus, by (\ref{ACOI1a}) holding on the absorbing set, the following will apply almost surely with $X_0=x$ where $x \in B$:
\begin{eqnarray}
&&E^{\gamma^*}[h(X_t) | x_{[0,t-1]},u_{[0,t-1]}] = \int_y h(y) P(X_t \in dy|x_{t-1},u_{t-1}) \\
&&\quad \quad \quad = c(x_{t-1},f(x_{t-1})) + \int_y h(y) P(dy|x_{t-1},f(x_{t-1})) - c(x_{t-1},f(x_{t-1})) \\
&& \quad \quad \quad \leq g + h(x_{t-1})  - c(x_{t-1},f(x_{t-1})) 
\end{eqnarray}
Iterating the above and dividing by $n$, we arrive at
\[g - {1 \over n} E_x^{\gamma^*}[h(X_n)] + {1 \over n} E_x^{\gamma^*}[h(X_0)] \geq {1 \over n}E_x^{\gamma^*}[\sum_{t=1}^n c(X_{t-1},U_{t-1})].\]
Taking the limsup on both sides (and replacing $\limsup$ with $\liminf$ by reversing the negative sign on the left), and (\ref{condConv0}) holding for $\gamma^* = \{f,f,f,\cdots\}$, we establish the desired bound.
\end{proof}

In particular if we have that $g$ is a lower bound on the optimal cost (say $g = \gamma^*$ in (\ref{E-thmWeak0}) as a consequence of the convex analytic method), we can claim that $\gamma^*$ is optimal for all initializations $X_0=x$ where $x \in B$.
 
Now, recall that the analysis in \cite[Theorem 5.3]{hernandez1993existence} shows that if we have an optimal invariant measure, then this leads to (\ref{ACOI11}) for some randomized $\phi$ on a set of measure 1 under $\pi_{\phi}$ with $h$ bounded from below. Building on \cite{Blackwell2,Blackwell3}, via \cite[(5.7)]{hernandez1993existence}, this implies the existence of a deterministic control policy $k$ which is defined on $B$ and which satisfies
\begin{align}\label{ACOI1}
g + h(x) \geq \bigg(c(x,k(x)) + \int h(x') {\cal T}(dx'|x,k(x)) \bigg) 
\end{align}
However, with $\kappa^* = \{k,k,k,\cdots\}$, to be able to claim the optimality of $k$ over $B$ via Lemma \ref{ACOIEqn}, we need to show $P^{\kappa^*}(x,B) = 1$ for all $x \in B$; that is, an absorbing set under $k$ should be a subset of the absorbing set under $\phi$ when $X_0 =x$ with $x \in B$.

If the induced Markov chain under $\phi$ is positive Harris recurrent, then \cite[Theorem 5.3(b)]{hernandez1993existence} shows that (\ref{ACOI11}) holds everywhere (that is, for all $x \in \mathbb{X}$), and the result follows. 

Additionally, when $\mathbb{U}$ is countable, this result also follows via the following argument: By Blackwell's theorem \cite[p. 864]{Blackwell2} and by the measurable selection theorem of Blackwell and Ryll-Nardzewski \cite{Blackwell3} (see also p. 255 of \cite{dynkin1979controlled}), $k$ can be (without loss) constructed such that for all $x: k(x) \in \{u: (c(x,u) +  \int h(x') {\cal T}(dx'|x,u)) \leq c(x,\phi(x)) + \int h(x') {\cal T}(dx'|x,\phi(x))  \} \cap \{u: \phi(u|x) > 0\}$. In this case, it follows by expressing the transition probabilities in terms of the countable collection of control realizations, we will have that $P^{\kappa^*}(x,B) = 1$ for all $x \in B$.



\begin{theorem}\label{optimalDeterministicStatHL}
\begin{itemize}
Assume that one of the following holds: \hyperlink{H1}{(H1)} holds and the bounded cost function $c(x,u)$ is continuous; or \hyperlink{H2}{(H2)(a)} holds and the bounded cost function $c(x,u)$ is continuous in $u$ for every $x$. Accordingly, either Theorem \ref{thmWeak} or Theorem \ref{thmSetwise} apply. Let $\mu_*$ be an optimal invariant measure.
Define a stationary policy $\gamma$ via the disintegration
\begin{equation}\label{denklem11defineP2}
 \mu_*(\D{x},\D{u})\,=\,\gamma_*(\D{u} \smid x)\,\uppi_*(\D{x})
 \end{equation}
$\mu_*$ almost surely. Take $\phi\in\Usm$ be any policy that agrees with $\gamma_*$ on the support of $\uppi_*$. 
\item[(i)]\cite{hernandez1993existence}  If the induced Markov chain under $\phi$ is positive Harris recurrent, then the optimal policy can be assumed to be deterministic. 
\item[(ii)] If the induced Markov chain under an optimal policy is not positive Harris recurrent, then with $\mathbb{U}$ countable, on the support of $\uppi_*$, $\phi$ can be assumed to be deterministic. This would lead to an optimal policy for all initial states $x$ with $X_0=x$ where $x \in \supp(\uppi_*)$.
\end{itemize}
\end{theorem}

\subsubsection{A fixed point theorem approach}

  
Borkar \cite{Borkar2} utilizes Schauder's fixed point theorem to arrive at the optimality of deterministic policies directly via the convex analytic method.


\begin{assumption}\label{A3.10}
There exists a $\sigma$-finite non-negative measure $\lambda$ on $\mathbb{X}$ such that
\begin{equation*}\cT(\D{y} \smid x,u) = f(x,u,y) \lambda(\D{y}),
\qquad x \in \XX\,,\ u \in \Act\,,
\end{equation*}
$f$ is continuous in all its variables, and $f(x, u, \cdot)$ is bounded and equicontinuous
(over $x \in \XX$, $u \in \Act$) and bounded away from zero
uniformly over all compact sets. The state and control variables are finite dimensional real valued. Furthermore ${\cal G}$ is compact and every stationary and randomized policy leads to a Markov chain which admits an invariant probability measure.
\end{assumption}

The conditions above are needed in order to employ a version of Schauder's fixed point
theorem on maps on the space of probability measures under the total variation distance. The above then leads to the following extremal property for deterministic policies.

\begin{theorem}[Lemma 11.16 in \cite{Borkar2}]
Under \cref{A3.10}, suppose that with $a \in (0,1)$ and $\phi^1, \phi^2$ two stationary control policies
\begin{equation*}
\phi(\D{u}\smid x) = a \phi^1(\D{u}\smid x) + (1-a) \phi^2(\D{u}\smid x)\,,
\end{equation*}
where $\phi^1(\D{u}\smid x) \neq \phi^2(\D{u}\smid x)$ for all $x$ (which can be refined to $x \in B_R$ for some ball of sufficiently large radius $R$). Then, the invariant probability measure induced by $\phi$ cannot be an extreme point.
\end{theorem}

\subsubsection{An approach via the small/petite set theory}\label{S3.5}

For completeness, we present an approach via the theory of small/petite sets to arrive at complementary conditions for the optimality of stationary and randomized policies in the appendix. The approach is to follow the proof method utilized in \cref{extremePoints} where the analysis reduces to a stochastic realization condition, which however does not appear to be lenient. The details are reported in the appendix: the realization condition itself is likely a useful property for further applications and for this reason the analysis is reported in the appendix.




\section{Denseness of Performance of Stationary Deterministic Policies}

In some applications it may be useful to know not only that optimal policies
are deterministic, but that deterministic policies are dense in the sense of
approximability of the costs induced under randomized and stationary policies. We will in fact show that the performance of deterministic, but also quantized (i.e., those with finite range), policies are dense. 



We have the following supporting \emph{denseness} result involving measurable policies over randomized ones.

\begin{theorem}\label{DenseQuantized} 
Let $(X,U)$ be finite dimensional real valued state and control action random variables,
where the compact $\Act$ valued $U$ is generated by a randomized stationary policy.
Suppose further that $X$ admits a non-atomic probability measure.
Then we have the following:
\begin{enumerate} 
\item[\ttup{i}] 
There exists a collection of measurable policies $U_n = \gamma_n(X)$ so that $(X,U_n)$ converges weakly to $(X,U)$.
\item[\ttup{ii}] $(X,U_n) \to (X,U)$ in the $w$-$s$ (setwise-weak) topology also.
\item[\ttup{iii}] If  $X_n$ is a sequence of random variables whose associated probability
measure converges in total variation to that of $X$, then the joint random variable
$(X_n,\gamma_n(X_n))$ converges weakly to $(X,U)$ as well as in the $w$-$s$ sense (setwise in $x$ and weakly in $u$).
\end{enumerate}
\end{theorem}

\begin{proof}
\begin{itemize}
\item[(i)] is due to \cite[Theorem 3]{milgrom1985distributional}, though there exist other related results, e.g. \cite[Proposition 2.2]{beiglbock2018denseness},
\cite{lacker2018probabilistic}, \cite[Theorem 3]{milgrom1985distributional},
but also many texts in optimal stochastic control where denseness of deterministic
controls have been established inside the set of relaxed controls \cite{borkar1988probabilistic}.
\item[(ii)] The marginal on $X$ is fixed along the sequence. The result then follows from \cite[Theorem 3.10]{schal1975dynamic} (or \cite[Theorem 2.5]{balder2001}).
\item[(iii)] Let $\rho$ denote the Prohorov metric on the joint
state-action random variables. Write
\begin{equation*}
\rho\bigl((X_n,\gamma_n(X_n)) , (X,U)\bigr)
\,\le\, \rho\bigl((X_n,\gamma_n(X_n)) , (X,\gamma_n(X))\bigr) + \rho\bigl((X,\gamma_n(X)) , (X,U)\bigr)\,.
\end{equation*}
The first term on the right converges to zero due to total variation convergence of $X_n$ to $X$ (since we apply the same deterministic measurable policy $\gamma_n$, and convergence is uniform over all measurable functions as in the proof of \cite[Lemma 1.1(iii)]{KYSICONPrior}). The second term converges to zero by (i). 

As in the proof of \cref{thmSetwise}, by \cite[Theorem 3.10]{schal1975dynamic} or \cite[Theorem 2.5]{balder2001}, since the measure converges weakly and the marginal in $X$ converges setwise, the convergence is also in the $w$-$s$ sense. 
\end{itemize}
\end{proof}

In fact, from the proof of \cref{DenseQuantized}(i) (see e.g.  \cite[Theorem 3]{milgrom1985distributional}) one shows not only the denseness of deterministic policies, but also those with a quantized range so that $| \gamma_n(\mathbb{X}) | < \infty$. 

\cref{DenseQuantized} helps us in establishing the following. 


\begin{theorem}\label{denseDet}
Suppose that 
\begin{itemize}
\item[\ttup{i}]
$\cG$ is weakly compact. Furthermore $\mathbb{X}=\mathbb{R}^n$ for some finite $n$, and for all $x \in \mathbb{R}$, $\mathbb{U}(x) = \mathbb{U}$ is compact.
\item[\ttup{ii}]
For some $\alpha \in [0,1)$, under every stationary policy
$\gamma$ the induced  kernel $P^{\gamma}$ of the Markov chain given by
\begin{equation*}
P^{\gamma}(\pi) (\cdot)\,\df\, (\pi {\cal T}^{\gamma})(\cdot) \,=\, \int \pi(\D{x})
\gamma(\D{u}\smid x) \int  \cT(\cdot \smid  x,u)
\end{equation*}
satisfies
\begin{align}\label{dobCont}
 \bnorm{P^{\gamma}(\pi) - P^{\gamma}(\bar{\pi})}_{\mathsf TV} \leq \alpha \|\pi-\bar{\pi}\|_{\mathsf TV}
 \end{align}
for any pair of probability measures $(\pi,\bar{\pi})$. This condition implies, naturally, that every stationary policy leads to a unique invariant probability measure.

\item[\ttup{iii}] The kernel $\cT(\D{y} \smid x,u)$ is such that, the family of conditional probability measures $\{\cT(\D{y} \smid x,u), x \in \mathbb{X}, u\in \mathbb{U}\}$ admit densities $f_{x,u}(y)$ with respect to a reference measure and all such densities are bounded and equicontinuous (over $x \in \mathbb{X}, u\in \mathbb{U}$). 

\item[\ttup{iv}] One of the following holds: \hyperlink{H1}{(H1)} holds and the bounded cost function $c(x,u)$ is continuous; or \hyperlink{H2}{(H2)(a)} holds and the bounded cost function $c(x,u)$ is continuous in $u$ for every $x$.
\end{itemize}
Then, deterministic and stationary policies are dense among those that are randomized and stationary, in the sense that the cost under any randomized stationary policy can be approximated arbitrarily well by deterministic and stationary policies. Furthermore, the dense set of deterministic and stationary policies can be assumed to have finite range. 
\end{theorem}

Before presenting the proof, we note that a list of sufficient conditions for \cref{dobCont}
are presented in \cite[Theorem 3.2]{hernandez1991recurrence} and these all have a
relationship with the Dobrushin's ergodicity coefficient \cite{dobrushin1956central}.

\begin{proof} Observe that the family of densities $f_{x,u}(\cdot)$ being equicontinuous over $x \in \mathbb{X}, u\in \mathbb{U}$ implies that $\{ \int f_{x,u}(y) \mu(dx)\gamma(du|x), \quad \mu \in {\cal P}(\mathbb{X}), \gamma \in \Gamma_{\mathsf S}\}$ is also equicontinuous. Thus, the family of invariant probability measures under any stationary policy admit densities (with respect to a reference measure) which are bounded and equicontinuous.
Then, following e.g., \cite[Lemma 4.3]{YukselOptimizationofChannels}, if $\{f_n\}$ is a sequence of probability density functions (with respect to some reference measure $\psi$) which are equicontinuous and uniformly bounded and if $\mu_n(dy) = f_n(y)\psi(dy) \to \mu(dy) = f(y) \psi(dy)$ weakly, then as a consequence of the Arzel\'a-Ascoli theorem (applied to $\sigma$-compact spaces) $f_n \to f$ pointwise and by Scheff\'e’s theorem, $\mu_n \to \mu$ in total variation. 

Let $\gamma$ be any randomized policy. Suppose that this policy gives rise to an invariant probability measure $\uppi_\gamma(\D{x},\D{u})$. Now, consider a sequence of deterministic policies $f_n$ so that under this sequence of policies $\uppi_\gamma(\D{x}) \delta_{f_n(x)}(\D{u})$ converges weakly to $\uppi_\gamma(\D{x},\D{u})$ by Theorem 4.1(i).

Now, let us apply the same measurable policy sequence to the random variable $X_n$ which has the probability measure $\uppi_{f_n}(\D{x})$ equal to the marginal of the invariant measure under policy $U=f_n(X)$. Then, for every continuous and bounded $g \in C_b(\XX)$
\begin{equation*}
\int \uppi_{f_{n}}(\D{x})\delta_{f_{n}(x)}(\D{u}) \bigg( \int g(y) \cT(\D{y} \smid  x,u) \bigg) = \int \uppi_{f_{n}}(\D{x}) g(x)\,.
\end{equation*}
Let $\uppi_{f_{n_k}}(\D{x})$ be a weakly converging subsequence with limit $\eta$ (by the compactness assumption on $\cG$). By hypothesis, this convergence is also in total variation. 


Define for any stationary policy $f$
\begin{equation*}
P^{f}(\pi) = \int \pi(\D{x})
f(\D{u}|x) \int  \bigg(\cT(\D{y} \smid  x,u) \bigg)\,,
\end{equation*}
and by hypothesis note that $\|P^{f}(\pi) - P^{f}(\bar{\pi})\|_{\mathsf TV} \leq \alpha \|\pi - \bar{\pi}\|$. Then, 
\begin{align*}
\|\uppi_{f_{n_k}} - \uppi_\gamma \|_{\mathsf TV} 
&= \| P^{f_{n_k}}(\uppi_{f_{n_k}}) - P^{\gamma}(\uppi_\gamma)\|_{\mathsf TV} \\
& = \| P^{f_{n_k}}(\uppi_{f_{n_k}}) - P^{f_{n_k}}(\uppi_\gamma)\|_{\mathsf TV} + \|  P^{f_{n_k}}(\uppi_\gamma) - P^{\gamma}(\uppi_\gamma) \|_{\mathsf TV} \\
&\leq \alpha \|\uppi_{f_{n_k}} - \uppi_\gamma \|_{\mathsf TV} + \|P^{f_{n_k}}(\uppi_\gamma) - P^{\gamma}(\uppi_\gamma)\|_{\mathsf TV}\,,
\end{align*}
and thus
\begin{equation*}\|\uppi_{f_{n_k}} - \uppi_\gamma \|_{\mathsf TV} \leq \frac{ \|P^{f_{n_k}}(\uppi_\gamma) - P^{\gamma}(\uppi_\gamma)\|_{\mathsf TV}}{1-\alpha} \end{equation*}
Now, for the right hand side, we have that $P^{f_{n_k}}(\uppi_\gamma) \to P^{\gamma}(\uppi_\gamma)$ weakly, since $\uppi_\gamma(\D{x})
f_n(\D{u}|x)$ converges in the $w$-$s$ sense to $\uppi_\gamma(\D{x}) \gamma(\D{u}|x)$ by \cref{DenseQuantized}(ii) and ${\cal T}$ is weakly continuous or strongly continuous in actions, under \ttup{iv}. 

However, by the discussion at the beginning of the proof above, this convergence also holds in total variation: Note that $P^{f_{n_k}}(\uppi_\gamma)$ is the measure defined only on the state marginal. This converges weakly but by the assumption of equicontinuity on the densities, the sequence of densities will have a converging subsequence. Therefore, the weak convergence should be supported by pointwise convergence of densities, and thus by Scheff\'e's lemma, the convergence is also in total variation.

Since the right hand side converges to zero, we can conclude that indeed $\uppi_{f_{n_k}}(\D{x}) \to \uppi_\gamma(\D{x})$. 

Finally, by \cref{DenseQuantized}(iii), as 
$\uppi_{f_{n_k}}(\D{x})\delta_{f_n(x)}(\D{u}) \to \uppi_\gamma(\D{x}) \gamma(\D{u}\smid x)$ in the setwise-weak (setwise in the state, weakly in the control action), the result follows. 

Since we can assume that the deterministic policies converging to the randomized policies in \cref{DenseQuantized}(i) have finite range, the proof is complete.

\end{proof}

We note that \cite[Theorem 3.2]{SaldiLinderYukselTAC14} and \cite[Theorem 4.2]{saldi2014near} had established near optimality of quantized policies (though not necessarily deterministic), under \hyperlink{H2}{(H2)} (with slightly more restrictive conditions) and  \hyperlink{H1}{(H1)}, respectively. The unified analysis here is more direct and general.

\section{Conclusion}
We presented results on the existence of optimal policies in average cost optimal stochastic control with kernels that do not satisfy weak kernel continuity in both state and actions, but with strong kernel continuity in the actions for every fixed state variable. We also studied conditions for the optimality of deterministic policies in average cost optimal stochastic control and reviewed some prior work. We finally presented a denseness result of costs induced under deterministic and stationary policies among those that are attained by randomized and stationary policies.



\appendix

\section{An approach based on the theory of small sets and an open realizability question}

In the following we present a sufficient condition to establish the desired optimality result on deterministic policies through the theory of small sets. 

%

\begin{definition}
A set $A \in \mathcal{B}(\XX)$ is n-small on
$\bigl(\XX,\mathcal{B}(\XX)\bigr)$ if for some positive measure $\mu_n$
$$P^n(x,B) \,\ge\, \mu_n(B), \quad \forall x \in A \text{\ and\ } B
\in \mathcal{B}(\XX)\,.$$
\end{definition}


\begin{definition} \cite{MeynBook}
A set $A \in \mathcal{B}(\XX)$ is $\nu_{\mathcal{K}}$-petite on
$\bigl(\XX,\mathcal{B}(\XX)\bigr)$, if for
some distribution $\mathcal{K}$ on $\mathbb{N}$ (set of natural numbers),
and some positive measure $\nu_{\mathcal{K}}$,
$$\sum_{n=0}^{\infty} P^n(x,B) \mathcal{K}(n) \,\ge\, \nu_{\mathcal{K}}(B)\,,
\quad \forall x \in A \text{\ and\ } B\in \mathcal{B}(\XX)\,.$$
\end{definition}

By  \cite[Proposition~5.5.6]{MeynBook}, if a Markov chain is $\psi$-irreducible,
 and if a set $C$ is $\nu$-petite, then $\mathcal{K}$ can be taken to be a geometric distribution $a_{\epsilon}(i) =(1 - \epsilon) \epsilon^{i}, \quad i \in \mathbb{N}$ (with the randomly sampled chain also known as the resolvent kernel).
 
{\bf The \texorpdfstring{$1$}{}-small case.} We impose the following assumption.

\begin{assumption}\label{SmallSetCond}
For any two policies $\gamma^1$ and $\gamma^2$
in $\Usm$ (possibly randomized), and every Borel set $B$ that satisfies
$\psi_{\gamma^i}(B)>0$, $i=1,2$, where $\psi_{\gamma^i}$ denotes
the maximal $\psi$-irreducibility measure under policy $\gamma^i$,
there exists a measurable $C \subset B$ that is a  $1$-small set with positive maximal irreducibility measure under
either of the transition probabilities $P^{\gamma^1}$ and $P^{\gamma^2}$. 
\end{assumption}

\begin{proposition}
A sufficient condition for \cref{SmallSetCond} is that the following hold:
\begin{enumerate}
\item[\ttup{i}] The transition kernel $\cT$ is bounded from below by a
conditional probability measure that admits a density with respect to some positive
measure $\phi$.
In other words there exist a measurable
$f\colon\XX\times\Act\times\XX\to\RR_+$,
such that
\begin{equation*}
\cT(D\smid x,u) \,\ge\, \int_D f(x,u,y) \phi(\D{y})
\end{equation*}
for every $D\in\mathcal{B}(\XX)$.
\item[\ttup{ii}] The function $f(x,u,y)$ in \ttup{i}
is continuous in $x, u$ for every fixed $y$, and $\Act$ is compact.
\item[\ttup{iii}]
It holds that
\begin{equation*}
\int_{\XX}\inf_{x \in A, u \in \Act}\,  f(x,u,y) \phi(\D{y})\,>\,0
\end{equation*}
for every nonempty compact set $A\subset\XX$.
\end{enumerate} 
\end{proposition}

\begin{proof}
The measurable selection results in
\cite{Schal,kuratowski1965general} and \cite[Theorem 2]{himmelberg1976optimal}
show that,
for any compact $A\subset\XX$, there exist
measurable functions $g$ and $F$ such that
\begin{equation*}
\inf_{x \in A, u \in \Act}f(x,u,y)
\,=\, \min_{x \in A, u \in \Act}f(x,u,y) \,=:\, F\bigl(g(y),y\bigr)
\end{equation*}
Thus, using the notation in \cref{E-Tgamma}, we have
\begin{equation*}
P^{\gamma}(x,D) \ge \int_{D} \inf_{x \in A, u \in \Act}\, f(x,u,y)  \phi(\D{y}) = \int_D F\bigl(g(y),y\bigr) \phi(\D{y}) =: \nu(D)
\end{equation*}
for some finite (sub-probability) measure $\nu$. Thus, every compact set is $1$-small under a given policy.
\end{proof}


\begin{theorem}
Under \cref{SmallSetCond}, and the realizability condition given in
\cref{targetD1,targetD2} (presented further below) a randomized policy cannot lead to an extreme measure in $\cG$, that is,
all measures in $\cG_e$ are induced by
deterministic policies. 
\end{theorem}

\begin{proof} 
The proof is divided into four steps.

\noindent{\bf Step 1.} Let there be a policy $\phi$ which is randomizing between two 
policies $\phi^1$ and $\phi^2$ on some measurable set $B$,
so that for some $\kappa(x) \in (0,1)$
with $x \in B$, we have that
\begin{equation*}
\phi(\D{u}\smid x) \,=\, \kappa(x) \phi^1(\D{u}\smid x)
+ (1- \kappa(x)) \phi^2(\D{u}\smid x)\,, \qquad x \in B\,.
\end{equation*}
By \cref{SmallSetCond}, there exists a $C \subset B$ so that this set is small
for either of the transition probabilities and on this set the above randomization
also holds, and that the measure on $C$ is positive under either of the irreducibility measures under $P^{\phi^1}$ and $P^{\phi^2}$.

We can assume that the transition kernels admit small sets with measure $\nu^1$ and $\nu^2$, where we take $\nu^1(\XX)=\nu^2(\XX)$, without any loss of generality, since we can always scale down the measure with the larger total mass to match the one with the smaller total mass.


\noindent{\bf Step 2.}
Define, for $K \in (0,\frac{1}{2})$,
\begin{equation*}
C^K\,\df\, \{x \in C:  1 - K \,\ge\, \kappa(x) \,\ge\, K \}\,.
\end{equation*}
Thus, we have
\begin{equation*}\phi(\D{u}\smid x) = \kappa(x) \phi^1(\D{u}\smid x) + (1- \kappa(x)) \phi^2(\D{u}\smid x), \qquad x \in C^K, \end{equation*}
and $C^K$ is also small (since it is a subset of a small set). Furthermore, we can take $C^K$ be so that it has positive measure under the irreducibility measures (by a continuity of measures argument, as $K \to 0$, $\psi(C^K) \to \psi(C)$ for any measure $\psi$). Now, by the Nummelin-Athreya-Ney split chain argument
\cite{Nummelin,NummelinSplit,AthreyaNey}, we split $C^K$ into $C^K \times \{0\}$ and $C^K \times \{1\}=:\alpha$, where $\alpha$
is a \emph{pseudo-atom}, in the sense that the transition kernels
are independent of the particular $x\in\alpha$, as we make more explicit below. To motivate this construction, we note that for $x \in C^K$, for any Borel $A$, we can write
\begin{align*}
\cT^{\phi}(A \smid x)  &\,=\, \biggl(\bigl(1-K\nu^1(\XX)\bigr)
\frac{\kappa(x)\cT^{\phi^1}(A\smid x) - K \nu^1(A)}{1-K\nu^1(\XX)}
 + K\nu^1(\XX) \frac{K \nu^1(A)}{K\nu^1(\XX)} \biggr)   \\ 
& \mspace{20mu}+   \biggl(\bigl(1-K\nu^2(\XX)\bigr)
\frac{ (1- \kappa(x))\cT^{\phi^2}(A\smid x) - K \nu^2(A)}{1-K\nu^2(\XX)} 
+ K\nu^2(\XX)\frac{K \nu^2(A)}{K\nu^2(\XX)}\biggr)\,.
\end{align*}
Write the above as
\begin{equation*}
\begin{aligned}
&\cT^{\phi}(A \smid x)\\
&= \frac{1}{2}\biggl(\bigl(1-2K\nu^1(\XX)\bigr)
\frac{ \Bigl(  \kappa(x)\cT^{\phi^1}(A\smid x) +  (1- \kappa(x))\cT^{\phi^2}(A\smid x)\Bigr) 
-  K \bigl(\nu^1(A) + \nu^2(A)\bigr)}{1-2K\nu^1(\XX)}  \\ 
& \mspace{480mu} + 2K\nu^1(\XX) \frac{2K \nu^1(A)}{2K\nu^1(\XX)} \biggr) \nonumber \\ 
& +  \frac{1}{2}\biggl(\bigl(1-2K\nu^2(\XX)\bigr)
\frac{ \Bigl(  \kappa(x)\cT^{\phi^1}(A\smid x) + (1- \kappa(x))\cT^{\phi^2}(A\smid x)\Bigr)
-  K \bigl(\nu^1(A) + \nu^2(A)\bigr)}{1-2K\nu^2(\XX)}   \\ 
& \mspace{480mu} + 2K\nu^2(\XX)\frac{2K \nu^2(A)}{2K\nu^2(\XX)} \biggr) \nonumber 
\end{aligned}
\end{equation*}

 
\noindent{\bf Step 3 (The realizability step).} 
We  now realize (i.e., construct) two control policies, called $\tilde{\phi}^1$ and $\tilde{\phi}^2$, so that these policies agree with $\phi$ outside $C^K$, and inside $C^K$ they
 admit a split chain where the transitions outside the atom $\alpha$, that is on $C^K \times \{0\}$, are also in agreement: the only difference is on the atom itself, therefore, the policies act as if they are in agreement everywhere except on the atom. That is, for $x \notin C^K$, we have
\begin{equation*}
\tilde{\phi}^1(\D{u}\smid x) \,=\, \tilde{\phi}^2(\D{u}\smid x)=\phi(\D{u}\smid x)
\end{equation*}

But on $C^K$, we have that $\cT^{\phi}(\,\cdot\,\smid x)$ is attained by randomizing between $\tilde{\phi}^1$ and $\tilde{\phi}^2$ according to:
\begin{equation}\label{splitRandomP}
\phi(\D{u}\smid x) = \frac{1}{2} \tilde{\phi}^1(\D{u}\smid x)
+ \frac{1}{2} \tilde{\phi}^2(\D{u}\smid x), \qquad x \in C^K
\end{equation}
where in the split chain, for $x \in C^K \times\{0\}$, $\tilde{\phi}^i(\D{u}\smid x)$ leads to the one step transition kernel
\begin{equation*}
\cT^{\tilde{\phi}^i}(\,\cdot\, \smid x) 
\,=\,  \kappa(x) \cT^{\phi^1}(\,\cdot\,\smid x) + (1- \kappa(x)) \cT^{\phi^2}(\,\cdot\,\smid x)
 - K \nu^1(\,\cdot\,) - K\nu^2(\,\cdot\,)
\end{equation*}
for $i=1,2$. And on $x \in \alpha$, $\tilde{\phi}^i(\D{u}\smid x)$ leads to the one step transition kernel
\begin{equation*}
\cT^{\tilde{\phi}^i}(\D{y}\smid x) \,=\,  \frac{\nu^i(\D{y})}{\nu^i(\XX)}
\end{equation*}
for $i=1,2$. Note that the above lead to the following virtual \emph{aggregate} transition kernels under $\tilde{\phi}^i$, $i=1,2$ for $x \in C^K$:
\begin{align}
\cT^{\tilde{\phi}^1}(\,\cdot\, \smid x)
&\,=\, \kappa(x)\cT^{\phi^1}(\,\cdot\,\smid x)
+ (1- \kappa(x)) \cT^{\phi^2}(\,\cdot\,\smid x)
+ K\nu^1(\,\cdot\,) - K\nu^2(\,\cdot\,)\,, \label{targetD1}  \\
\cT^{\tilde{\phi}^2}(\,\cdot\, \smid x) 
&\,=\, \kappa(x) \cT^{\phi^1}(\,\cdot\,\smid x) + (1- \kappa(x))\cT^{\phi^2}(\,\cdot\,\smid x)
-  K\nu^1(\,\cdot\,) +  K \nu^2(\,\cdot\,) \label{targetD2}
\end{align}
so that \cref{splitRandomP} holds. 

\noindent{\bf The Realizability Condition:} {\it There exist stationary control policies $\psi^1$ and $\psi^2$, such that, $\psi^1$ realizes \cref{targetD1} and $\psi^2$ realizes \cref{targetD2}.}

As a result, the only difference in the transition kernels, as seen from the split chain/atom is that, in the atom $\alpha$ randomization occurs; outside the atom the transition probabilities are identical.

\noindent{\bf Step 4.}
$\alpha$ is the \emph{accessible atom} of interest: $\Exp_{\alpha}^{\tilde{\phi}^i}[\tau_{\alpha}]< \infty$ and $P^{\tilde{\phi}^i}(x,B)=P^{\tilde{\phi}^i}(y,B)$ for all $x,y \in \alpha$, where $\tau_{\alpha}=\min(k>0: x_k=\alpha)$ is the return time to $\alpha$. Since $\phi$ is randomizing between the two policies $\phi^1$ and $\phi^2$, let $v, v^1, v^2$ be corresponding invariant measures to $\phi$, $\tilde{\phi}^1$ (only different at the atom), $\tilde{\phi}^2$ (only different at the atom); as noted above, for states other than those in the atom the transition kernels are identical. 

In this case, the invariant measures computed through the mean
empirical occupation measures normalized with the expected return times is obtained through the following analysis: For every Borel $A \in \XX, B \in \Act$,
\begin{align*}
v(A,B) &\,=\, \frac{\Exp_{\alpha}^{\phi}[\sum_{k=0}^{\tau_{\alpha}-1} \Ind_{A\times B}(X_k,U_k)]}
{\Exp_{\alpha}^{\phi}[\tau_{\alpha}]}  \\
&\,=\,  \frac{\frac{1}{2}\Exp_{\alpha}^{\tilde{\phi}^1}
\bigl[\sum_{k=0}^{\tau_{\alpha}-1} \Ind_{A\times B}(X_k,U_k)\bigr]
+  \frac{1}{2}\Exp_{\alpha}^{\tilde{\phi}^2}
\bigl[\sum_{k=0}^{\tau_{\alpha}-1} \Ind_{A\times B}(X_k,U_k)\bigl]}
{\frac{1}{2}\Exp_{\alpha}^{\tilde{\phi}^1}[\tau_{\alpha}] 
+  \frac{1}{2}\Exp_{\alpha}^{\tilde{\phi}^2}[\tau_{\alpha}]} \\
&\,=\, \frac{\Exp_{\alpha}^{\tilde{\phi}^1}
\bigl[\sum_{k=0}^{\tau_{\alpha}-1} \Ind_{A\times B}(X_k,U_k)\bigr]}
{\Exp_{\alpha}^{\tilde{\phi}^1}[\tau_{\alpha}]  +  \Exp_{\alpha}^{\tilde{\phi}^2}[\tau_{\alpha}]}
+ \frac{\Exp_{\alpha}^{\tilde{\phi}^2}
\bigl[\sum_{k=0}^{\tau_{\alpha}-1} \Ind_{A\times B}(X_k,U_k)\bigr]}
{\Exp_{\alpha}^{\tilde{\phi}^1}[\tau_{\alpha}] + \Exp_{\alpha}^{\tilde{\phi}^2}[\tau_{\alpha}]} \\
&\,=\, \frac{{\Exp_{\alpha}^{\tilde{\phi}^1}[\tau_{\alpha}]}
{\Exp_{\alpha}^{\tilde{\phi}^1}[\tau_{\alpha}]  + \Exp_{\alpha}^{\tilde{\phi}^2}[\tau_{\alpha}]}}\,\frac{\Exp_{\alpha}^{\tilde{\phi}^1}
\bigl[\sum_{k=0}^{\tau_{\alpha}-1}\Ind_{A\times B}(X_k,U_k)\bigr]}
{\Exp_{\alpha}^{\tilde{\phi}^1}[\tau_{\alpha}]}
 \\
&\mspace{50mu}+
\frac{\Exp_{\alpha}^{\tilde{\phi}^2}[\tau_{\alpha}]}
{\Exp_{\alpha}^{\tilde{\phi}^1}[\tau_{\alpha}]  +  \Exp_{\alpha}^{\tilde{\phi}^2}[\tau_{\alpha}]}\,
\frac{ \Exp_{\alpha}^{\tilde{\phi}^2}
\bigl[\sum_{k=0}^{\tau_{\alpha}-1} \Ind_{A\times B}(X_k,U_k)\bigr]}
{\Exp_{\alpha}^{\tilde{\phi}^2}[\tau_{\alpha}]}\\
&\,=\, \frac{ \Exp_{\alpha}^{\tilde{\phi}^1}[\tau_{\alpha}]}
{\Exp_{\alpha}^{\tilde{\phi}^1}[\tau_{\alpha}]
+  \Exp_{\alpha}^{\tilde{\phi}^2}[\tau_{\alpha}] } \,v^1(A,B) 
 + \frac{\Exp_{\alpha}^{\tilde{\phi}^2}[\tau_{\alpha}]}
{\Exp_{\alpha}^{\tilde{\phi}^1}[\tau_{\alpha}]
+  \Exp_{\alpha}^{\tilde{\phi}^2}[\tau_{\alpha}] }\,v^2(A,B)\,,
\end{align*}
which is a convex combination of $v^1$ and $v^2$.
In the above, the second equality is
 critical for the validity of the convex combination:
$\Exp_{\alpha}^{\phi}[\tau_{\alpha}] = \frac{1}{2}\Exp_{\alpha}^{\tilde{\phi}^1}[\tau_{\alpha}] 
+  \frac{1}{2}\Exp_{\alpha}^{\tilde{\phi}^2}[\tau_{\alpha}]$,
which holds due to the construction in {\bf Step 3}.
\end{proof}

We note also that the similar program applies for a construction building on both
$m$-small sets and petite sets. These sets exist under much less stringent
conditions than those required on $1$-small sets.

\end{document}